\documentclass[acmsmall,screen]{MapleTrans}

\usepackage{etoolbox}
\usepackage{breqn}
\usepackage{lineno}

\BeforeBeginEnvironment{dmath}{\begin{nolinenumbers}}%
\AfterEndEnvironment{dmath}{\end{nolinenumbers}}

\usepackage{amsthm}
\usepackage{amsmath}
\usepackage{amsfonts,latexsym}

\usepackage{graphicx}
\usepackage{multicol}
\usepackage{color}
\usepackage{pgfplots}
\usepackage{tikz, tikz-3dplot, float}

\usepackage{algorithm}
\usepackage{algorithmicx}
\usepackage{algpseudocode}
\usepackage{verbatim}

\usepackage{lhelp}
\usepackage{caption}
\usepackage{cleveref}
\usetikzlibrary{tikzmark}
\usetikzlibrary{positioning,calc}
\usetikzlibrary{arrows.meta,chains,decorations.pathreplacing,scopes,positioning}
\usetikzlibrary{calc,shapes.multipart,chains,arrows}
\DeclareGraphicsExtensions{.eps,.png}
\usepackage{multirow}

\def\L {\ensuremath{\mathfrak{L}}}

\def\Q {\ensuremath{\mathbb{Q}}}

\renewcommand{\underline}[1]{{#1}}

\newcommand{\Maple}{{\sc  Maple}}

\newcommand{\hidetext}[1]{\mbox{ \ }}

\newcommand{\textblockcomment}[1]{}

\begin{document}
\title{A New Algorithm for Computing Integer Hulls of 2D Polyhedral Sets}

\author{Chirantan Mukherjee}
\affiliation{%
  \institution{The University of Western Ontario}
  \streetaddress{1151 Richmond St}
  \city{London}
  \country{Canada}}
\email{cmukher@uwo.ca}
\orcid{0000-0002-6692-3305}


 \begin{abstract}
 {\bf Abstract.} The \texttt{IntegerHull} function is part of Maple's \texttt{PolyhedralSets} library, which calculates the integer hull of a given polyhedral set. This algorithm works by translating the supporting hyperplanes of the facets of the input polyhedral set inwards till each hyperplane encounters at least one integer point. The polyhedral set is then divided into smaller regions and a brute force method is applied to find the remaining vertices.

There are certain edge case scenarios where the computational cost of the existing algorithm can be high, for which we propose a new algorithm. We translate the supporting hyperplanes of the facets inwards from the (relative) opposite vertex till the hyperplanes encounter at least one integer point. Then, we can follow the same procedure as the old algorithm or use a recursive technique on the smaller regions.

The edge case scenarios mentioned above occurs when there are integer points present on the supporting hyperplanes of the facets of the polyhedral set. This increases the region on which the brute force method is applied to find the remaining vertices.
 \end{abstract}


\maketitle

\section{Introduction}
\label{sec1}
The integer points of rational polyhedral sets are of great interest in various areas of scientific computing. Two such areas are {\verb|combinatorial optimization|} (in particular integer linear programming) and {\verb|compiler optimization|} (in particular, the analysis, transformation and scheduling of for-loop nests in computer programs), where a variety of algorithms solve questions related to the points with integer coordinates belonging to a given polyhedron.

Wang and Moreno Maza have developed the existing integer hull algorithm for polyhedral sets in Maple \cite{wang2022, 10.1007/978-3-030-85165-1_15}. Their algorithm is applicable to polyhedral sets of any given dimension. However, for the purpose of our discussion, we will focus on the 2D case.

In this paper we propose a novel algorithm~\ref{alg:NewIntegerHull} for computing the integer hull of 2D polyhedral sets. The existing integer hull algorithm, as referenced in \cite{wang2022, 10.1007/978-3-030-85165-1_15}, relies on a recursive construction, effectively reducing the computation of integer hulls in arbitrary dimensions to that of dimension $2$. However, to enhance performance, we propose improvements for this base case. Our preliminary experiments with the new algorithm demonstrate its efficiency and ability to handle polyhedral sets with a large number of integer point.

There are three main families of integer hull algorithms for polyhedral sets: {\verb|cutting plane method|}, {\verb|branch-and-bound method|} and {\verb|lattice point counting method|}.

The {\verb|cutting plane method|} was introduced by Gomory \cite{Gomory2010} to solve integer linear programming and mixed integer linear programming. This method involves solving the linear programming problem to find the optimal solution. It does this by introducing new constraint at each step until an integer solution is found, which in turn reduces the area of the feasible solution. Chv\'atal \cite{CHVATAL1973305} and Schrijver \cite{SCHRIJVER1980291} provided a geometric description of this method and developed a procedure for computing the integer hull of a polyhedral set.

Land and Doig \cite{Land1960AnAM} introduced the {\verb|branch-and-bound method|} for computing the integer hull of a polyhedral set. This method works by recursively partitioning the polyhedral set into sub-polyhedral sets, then computing the integer vertices of each of the sub-polyhedral set and finally merging them all together.

Pick's theorem \cite{pick1899} can be used for calculating the area of any polygon with integer lattice points. Using this idea Barvinok \cite{Barvinok1994} created an algorithm for counting the number of integer lattice points inside a polyhedron. Building on Barvinok's algorithm, Verdoolaege, Seghir, Beyls, Loechner and Bruynooghe \cite{Verdoolaege2007CountingIP} came up with a method for counting the number of integer points in a non-parametric polytope. Meanwhile, Seghir, Loechner and Meister \cite{Seghir2012} developed a method of counting the number of images by an affine integer transformation of the lattice points contained in a parametric polytope. In 2004, a software package named {\verb|LATTE|}  \cite{DELOERA20041273} was developed for counting number of integer points in a rational polyhedral set using Barvinok's algorithm.

The paper is organized as follows. Section~\ref{sec2}  and \ref{sec3} is a brief review of polyhedral geometry and the existing integer hull algorithm in Maple. Section \ref{sec4} and \ref{sec5} presents the new integer hull algorithm in Maple for 2D cases. Section \ref{sec6} reports on our preliminary experimentation with the proposed algorithm.

\section{Preliminaries} 
\label{sec2}


We denote by $\mathbb{Z}$, $\mathbb{Q}$ and $\mathbb{R}$ the ring of integers, the
field of rational numbers and the field of real numbers. Unless specified otherwise,
all matrices and vectors have their coefficients in $\mathbb{Z}$.

\begin{definition}
A {\verb|polyhedral set|} $P$ is a set which can be expressed as the intersection of finite number of (closed) half spaces, that is $\{{\bf x} \in\mathbb{R}^n\mid A{\bf x}\leq {\bf b}\}$, where $A\in\mathbb{R}^{m\times n}$ is a matrix and ${\bf b}\in\mathbb{R}^m$ is a vector. 
\end{definition}

\begin{definition}
A subset $F$ of the polyhedron $P$ is called a {\verb|face|} of $P$ if $F$ equals $\{{\bf x}\in P \mid A_{\text{sub}}{\bf x} = {\bf b}_{\text{sub}} \}$ for a sub-matrix $A_{\text{sub}}$ of $A$ and a sub-vector ${\bf b}_{\text{sub}}$ of ${\bf b}$. 
A face of $P$, distinct from $P$ and of maximum dimension, is called a {\verb|facet|} of $P$. A face of dimension $0$ is called a {\verb|vertex|} of $P$.
\end{definition}

\begin{definition}
Given a polyhedral set $P$ and a vertex $v$ of $P$, the {\em vertex cone}
of $P$ at $v$ is the intersection of the half-spaces defining
$P$ and whose boundaries intersect at $v$.
\end{definition}

In our study of two-dimensional polyhedral sets, we observe that every non-trivial face falls into one of two categories: either it is a {\verb|facet|} (which can be a segment or a half-line), or it is a {\verb|vertex|}. Furthermore, each vertex cone exhibits a simplicial structure. Specifically, a vertex cone is defined by two half-lines originating from the same point, that is the vertex of that particular cone.

\section{Existing Integer Hull Algorithm}
\label{sec3}

The integer points of polyhedral sets hold paramount importance in the context of the
delinearization problem \cite{Benabderrahmane2010ThePM, Feautrier1988}, scheduling for-loop nests \cite{Feautrier1996}, accessing memory location of for-loop nests \cite{koppe2007}, Barvinok's algorithm for counting integer points in a polyhedral set \cite{Barvinok1994} and many more.

\begin{definition}
  The {\verb|integer hull|} $P_I$ of a convex polyhedral set $P$ is the convex hull of integer points of $P$.
\end{definition}

In this section we focus on implementing the current integer hull algorithm in Maple as described in \cite{wang2022, 10.1007/978-3-030-85165-1_15}. It has three main steps:
\begin{itemize}
\item {\verb|Normalization|}: In this step, the polyhedron $P$ is transformed into a rational polyhedron  $Q \subseteq {\Q}^d$. This transformation ensures that each facet of $Q$  has integer points on its supporting hyperplane, while maintaining $P_I = Q_I$.
\item {\verb|Partitioning|}: Here, integer points within $Q$ are identified and used to partition $Q$ into smaller polyhedral sets. Each set's integer hull can then be computed more straightforwardly.
\item {\verb|Merging|}: This step involves merging the integer hulls obtained from the partitioning process using a convex hull algorithm.
\end{itemize}

\begin{figure}[H]
  \centering
    \begin{tikzpicture}[scale=0.7]
    \draw[help lines, color=gray!40, dashed] (-3.9,-1.9) grid (4.9,5.9);
    \draw[red!30, ultra thick] (-2,-0.2) -- (3,-0.2) -- (1.7,3.9) -- cycle;
    \draw[color=red!60, fill=red!5] (-1,0) -- (2,0) -- (2,2) -- (1,3) -- (0,2) -- cycle;
    \draw[black!30, ultra thick, <->] (-2.5,0) -- (3.4, 0);
    \draw[black!30, ultra thick, <->] (1.7-0.3,3.9) -- (3-0.3,-0.3);
    \filldraw[red] (-2,-0.2)  circle (2pt) node[anchor=east] {$A$};
    \filldraw[red] (3,-0.2) circle (2pt) node[anchor=west] {$B$};
    \filldraw[red] (1.7,3.9) circle (2pt) node[anchor=south] {$C$};
    \draw[red] (-1,0) circle (2pt) node[anchor=east] {$D$};
    \draw[red] (0,0) circle (2pt);
    \draw[red] (1,0) circle (2pt);
    \draw[red] (2,0) circle (2pt) node[anchor=west] {$E$};
    \draw[red] (2,1) circle (2pt);
    \draw[red] (2,2) circle (2pt) node[anchor=west] {$F$};
    \draw[red] (0,1) circle (2pt);
    \draw[red] (0,2) circle (2pt) node[anchor=east] {$H$};
    \draw[red] (1,1) circle (2pt);
    \draw[red] (1,3) circle (2pt) node[anchor=east] {$G$};
    \draw[red] (1,2) circle (2pt);
    \draw[black] (2.4, -0.2) node[anchor=north] {$B'$};
    \draw[black] (1.5, 3.7) node[anchor=east] {$C'$};
    \draw[black] (2.6, 0) node[anchor=south] {$B''$};
    \draw[black] (-2, 0) node[anchor=south] {$A'$};
    \end{tikzpicture}
    \captionof{figure}{Integer hull of a triangle.}\label{fig5}
    \end{figure}
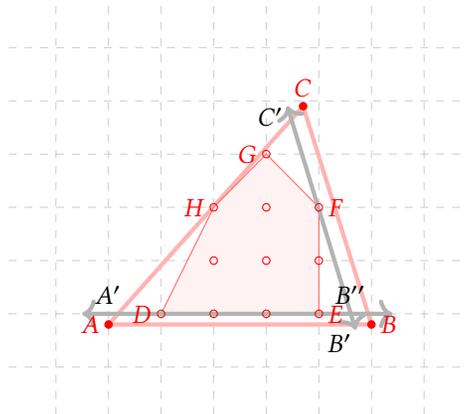

\begin{example}\label{ex1}
  The integer hull of the following triangle can be represented by the shaded pentagonal region in figure \ref{fig5}.
\end{example}

The integer hull $P_I$ of a convex polyhedral set $P$ can be described in terms
of its vertices. Since $P_I$ is the smallest polyhedral set containing
the integer points of $P$, the vertices of $P_I$ are necessarily integer
points.

For $P=\Delta ABC$ given in the above example \ref{ex1}, we have $\text{VertexSet}(P) = \{A, B, C\}$,
none of which are integer points. We translate the supporting hyperplane of the
facet $BC$ to the west until the supporting hyperplane of the facet $BC$
has at least one integer point, call it $F$ . We obtain a new triangle
$Q = \Delta AB'C'$, which, clearly has the same integer points as $P$.

We then repeat the process for the supporting hyperplane of the facet $AB'$ by translating
it upwards until we get integer points $\{D, E\}$, obtaining a
new facet $A'B''$.  We do not need to translate the supporting hyperplane of the
facet $A'C'$ since there is already an integer point $H$ on it.

We see that the triangle $\Delta A'B''C'$ can be divided into $4$ regions:
\begin{enumerate}
  \item the convex hull, say $R$,  of the points $\{D, E, F, H\}$;
  \item three "small" triangles:  $\Delta A'DH, \Delta HC'F, \Delta EB''F$.
\end{enumerate}
Since the vertices of $R$ are all integer points, $R$ is its own
integer hull. However, $R$ may not be the integer hull of $Q$. Indeed,
each of the small three triangles may still contain integer
points. This is actually the case for $\Delta HC'F$.

Because these three triangles are generally small in practice,
one can apply a brute force method to search for
integer points. One brute force method for computing the
integer hull $P_I$ of a polyhedral set $P$ is  to use
{\verb|Fourier-Motzkin elimination|} \cite{Banerjee1993,Pugh1992} to obtain
a parametric representation of $P$. With such a parametric
representation, one can enumerate all the integer points of $P$.

We want to stress that the computational cost of any brute force algorithm is inherently high, as such it is favorable to apply such a method only if the area is significantly small. In the next section we will discuss our alternative algorithm \ref{alg:NewIntegerHull} which is designed to ensure that the area on which the brute force method is applied remains small, thereby optimizing the computational efficiency.

Returning to our example, we complete the construction of the integer hull by putting together:
\begin{enumerate}
  \item the integer points found in each small triangle;
  \item the vertices of $R$;
  \item apply a convex hull algorithm, such as {\verb|QuickHull|} \cite{quickhull}, to all those points.
\end{enumerate}

\par

Let us first state the condition which guarantee the existence of integer points on any line \cite{wang2022, 10.1007/978-3-030-85165-1_15}.

A (parametric) polyhedral set can be defined as the intersection of vertex cones. 
That is a polyhedral set $P(\mathbf{b})$ defined by $\begin{pmatrix}
  a_1 & c_1\\
  a_2 & c_2\\
  \vdots & \vdots \\
  a_d & c_d
  \end{pmatrix} \begin{pmatrix}
    x\\
    y
    \end{pmatrix} \le \begin{pmatrix}
      b_1\\
      b_2 \\
      \vdots \\
      b_d
      \end{pmatrix}$ can be written as intersections of all vertex cones $S_i$, that is $P_i(\mathbf{b})=\bigcap\limits_{i=1}^{d} S_{i}(\mathbf{b})$, where each polyhedral cone $S_i$ is given by $\begin{pmatrix}
        a_i & c_i\\
        a_{i+1} & c_{i+1}\\
        \end{pmatrix} \begin{pmatrix}
          x\\
          y
          \end{pmatrix} \le \begin{pmatrix}
            b_i\\
            b_{i+1}\\
            \end{pmatrix}$, for all $i\in[1,d]$ such that $d+1 = 1$. 

The following lemma guarantees the existence of an integer vertex on the supporting hyperplane of the facets in the polyhedral set.

\begin{lemma}[\cite{10.1007/978-3-030-85165-1_15} Lemma 1 on p. 256, \cite{wang2022} Lemma 1 on p. 13]\label{lemma1}
  A line $ax+cy=b$, where $a,b,c\in\mathbb{Z}$ and $gcd(a,b,c)=1$, can have integer vertex $(x,y)$ if one of the following conditions is satisfied:
  \begin{enumerate}
    \item If $a\ne 0$ and $c\ne 0$, then $(x,y)$ is a integer vertex if and only if $gcd(a,c)=1$, that is, $x\equiv \frac{b}{a} (mod c)$.
    \item If $a=0$ (similarly $c=0$), then $(x,y)$ is a integer vertex if and only if $c=1$ (similarly $a=1$). 
  \end{enumerate}
\end{lemma}

Another crucial observation for constructing the integer hull is the ability to {\verb|translate|} the supporting hyperplane of each facet. This translation ensures that we can find at least one integer point on it. In other words, the integer hull of the parametric polyhedral set $P_I(\mathbf{b})$ remains unchanged when translating the supporting hyperplane of the facets by an integer $T$, that is, $P_I(\mathbf{b}+T)=P_I(\mathbf{b})$.

\par


\par



We can combine everything that we have observed so far in this section and construct the integer hull of any parametric 2D polyhedral set as follows.


We translate each supporting hyperplane $H$ of a facet $F$ of the polyhedral set $P$ inwards until $H$ intersects an integer point of the integer hull of the polyhedral set $P_I$. Indeed,
\begin{enumerate}
    \item We can detect whether $H$ has integer points or not by means of Lemma \ref{lemma1}. 
    \item If $H$ has integer points, Lemma \ref{lemma1} gives a formula for them, which we can then plug into the system of inequalities defining $P$, so as to check whether some of those integer points of $H$ are in $P_I$ or not.
\end{enumerate}

\section{New Integer Hull Algorithm}
\label{sec4}

In this section, we present our algorithm for computing the integer hull of a 2D polyhedral set. The integer hull algorithm proposed by Wang and Moreno Maza is extremely efficient as evidenced by the data presented in [\cite{wang2022} Table 4.5 on p. 51, Table 4.6 on p. 52]. However, it is crucial to note that when there are integer points present on the supporting hyperplane of the facets of the polyhedral set the computational cost become high. This observation becomes evident when considering the following example.

\begin{example}\label{limitations}
Consider the following triangle $P=\Delta ABC$, where none of whose vertices are integer points.
Following the algorithm discussed in the last section, we translate the supporting hyperplane of the facet $AB$ upwards until it intersects at least one integer point, call it $D$. This leads to a new triangle $\Delta A'B'C$, which retains the same integer points as $P$. We repeat this process for the supporting hyperplane of the facet $A'C$ by translating it to the east until we obtain integer points $\{D,H\}$ and a new triangle $Q=\Delta A''B'C'$. No adjustment is needed for the supporting hyperplane of facet $B'C'$ since it already contains an integer point, $F$.
    
The resulting triangle $\Delta A''B'C'$ can be partitioned into $4$ regions:
  \begin{itemize}
    \item item the convex hull, denoted as $R$, of the points $\{D,H,F\}$
    \item three triangles $DHB'$, $HC'F$ and $FDB'$.
  \end{itemize}
While the vertices of $R$ are all integer points, it may not represent the integer hull of $Q$. 
Each of the smaller triangles, $DHB'$, $HC'F$ and $FDB'$ , may still contain additional integer points.
Therefore, to ensure completeness, a brute force procedure is employed to identify integer points within these triangles. 
However, it is noteworthy that given the larger area of these triangles the brute force procedure will be computationally expensive. 
The time for obtaining the number of points using the brute force approach will surpass those obtained using the algorithm. In fact, the situation gets worse if there are integer points in the other edges as well.

\end{example}

\begin{figure}[H]
  \centering
    \begin{tikzpicture}[scale=0.7]
    \draw[help lines, color=gray!40, dashed] (-3.9,-1.9) grid (4.9,5.9);

    \draw[red!30, ultra thick] (-2.5,-0.2) -- (2.2,-0.7) -- (3.6,3.9) -- cycle;
    \draw[color=red!60, fill=red!15] (-2,0) -- (1,2) -- (2.3,-0.4) -- cycle;
    \draw[color=red!60, fill=red!15] (3.6,3.8) -- (1,2) -- (3,2) -- cycle;
    \draw[color=red!60, fill=red!15] (-2,0) -- (2.3,-0.4) -- (3,2) -- cycle;

    \filldraw[red] (-2.5,-0.2)  circle (2pt) node[anchor=east] {$A$};
    \filldraw[red] (2.2,-0.7) circle (2pt) node[anchor=west] {$B$};
    \filldraw[red] (3.6,3.9) circle (2pt) node[anchor=south] {$C$};

    \draw[red] (-2,0) circle (2pt) node[anchor=south] {$D$};
    \draw[red] (-1,0) circle (2pt);
    \draw[red] (0,0) circle (2pt);
    \draw[red] (1,0) circle (2pt);
    \draw[red] (2,0) circle (2pt);
    \draw[red] (2,1) circle (2pt);
    \draw[red] (2,2) circle (2pt);
    \draw[red] (0,1) circle (2pt);
    \draw[red] (1,1) circle (2pt);
    \draw[red] (1,2) circle (2pt) node[anchor=east] {$H$};
    \draw[red] (3,3) circle (2pt);
    \draw[red] (3,2) circle (2pt) node[anchor=west] {$F$};
    \draw[black!30, ultra thick, <->] (-3,-0.2+0.25) -- (2.2+0.3,-0.7+0.25);
    \draw[black!30, ultra thick, <->] (-2.3,-0.2) -- (3.8,+3.9);
    \draw[black] (-2.1, 0.25) node[anchor=east] {$A'$};
    \draw[black] (3, -0.5) node[anchor=east] {$B'$};
    \draw[black] (-2, 0) node[anchor=west] {$A''$};
    \draw[black] (3.7, 3.8) node[anchor=north] {$C'$};

    \end{tikzpicture}
    \captionof{figure}{Wang and Moreno Maza's integer hull algorithm.}\label{fig6}
    \end{figure}
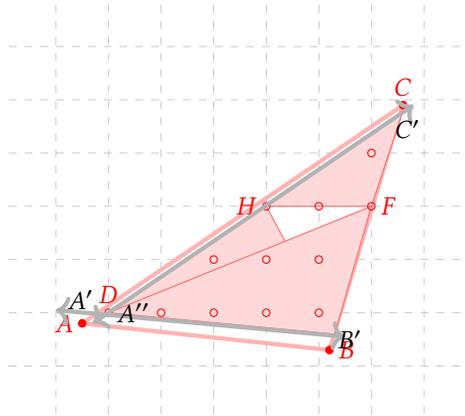

We propose a new algorithm to reduce the area on which the brute force method is applied.

We initiate the algorithm by translating the supporting hyperplane of facet $AB$ downward from the vertex $C$ until it intersects at least one integer point, denoted as $G$. This process is then repeated for the supporting hyperplanes of facets $AC$ and $BC$ by translating them to the west from vertex $B$ and east from vertex $A$, respectively, until integer points $E$ and $D$ are obtained.
  
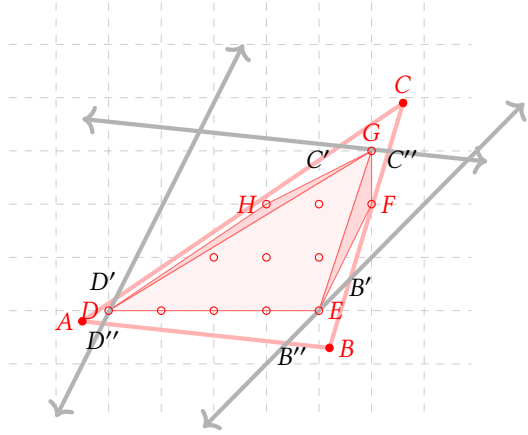
\begin{figure}[H]
\centering
  \begin{tikzpicture}[scale=0.7]
  \draw[help lines, color=gray!40, dashed] (-3.9,-1.9) grid (4.9,5.9);
      \draw[red!30, ultra thick] (-2.5,-0.2) -- (2.2,-0.7) -- (3.6,3.9) -- cycle;
      \draw[color=red!60, fill=red!5] (-2,0) -- (2,0) -- (3,3) -- cycle;
      \draw[color=red!60, fill=red!15] (-2,0) -- (1,2) -- (3,3) -- cycle;
      \draw[color=red!60, fill=red!15] (2,0) -- (3,2) -- (3,3) -- cycle;
      \draw[black!30, ultra thick, <->] (2.2-5.2,-2) -- (4-3.44,5);
      \draw[black!30, ultra thick, <->] (-2.5,-0.2+3.8) -- (2.2+3,-0.7+3.5);
      \draw[black!30, ultra thick, <->] (-2.5+2.3,-0.2-2) -- (3.6+2.3,+3.9);
      \filldraw[red] (-2.5,-0.2)  circle (2pt) node[anchor=east] {$A$};
        \filldraw[red] (2.2,-0.7) circle (2pt) node[anchor=west] {$B$};
        \filldraw[red] (3.6,3.9) circle (2pt) node[anchor=south] {$C$};
        \draw[red] (-2,0) circle (2pt) node[anchor=east] {$D$};
        \draw[red] (-1,0) circle (2pt);
        \draw[red] (0,0) circle (2pt);
        \draw[red] (1,0) circle (2pt);
        \draw[red] (2,0) circle (2pt) node[anchor=west] {$E$};
        \draw[red] (2,1) circle (2pt);
        \draw[red] (2,2) circle (2pt);
        \draw[red] (0,1) circle (2pt);
        \draw[red] (1,1) circle (2pt);
        \draw[red] (1,2) circle (2pt) node[anchor=east] {$H$};
        \draw[red] (3,3) circle (2pt) node[anchor=south] {$G$};
        \draw[red] (3,2) circle (2pt) node[anchor=west] {$F$};
        \draw[black] (-2.1, -0.2) node[anchor=north] {$D''$};
        \draw[black] (-2.1, 0.2) node[anchor=south] {$D'$};
        \draw[black] (1.5, -0.5) node[anchor=north] {$B''$};
        \draw[black] (2.8, 0.8) node[anchor=north] {$B'$};
        \draw[black] (2, 3.2) node[anchor=north] {$C'$};
        \draw[black] (3.6, 3.2) node[anchor=north] {$C''$};
      \end{tikzpicture}
  \captionof{figure}{New integer hull algorithm.}\label{fig7}
  \end{figure}

  The hexagon $D''B''B'C''C'D'$ resulting from this algorithm can be subdivided into $4$ regions:
  \begin{enumerate}
  \item the convex hull, denoted as $R$, formed by the points $\{D,E,G\}$
  \item three quadrilaterals $DD''B''E$, $EB'C''G$ and $D'DGC'$.
\end{enumerate}
While the vertices of $R$ are all integer points, $R$ may not represent the integer hull of $Q$. 
Indeed, each of the "small" quadrilaterals may still contain integer points, as is the case for $EB'C''G$ and $D'DGC'$.

To ensure that no integer points are missed, a brute force method can now be applied to search for integer points since the quadrilaterals are small. 
Alternatively, the algorithm can be repeated for quadrilaterals $EB'C''G$ and $D'DGC'$ until integer points $H$ and $F$ are obtained, adding these points to $\text{VertexSet}(P_I)$. 
For quadrilateral $DD''B''E$, there are no integer points.

Hence, $\text{VertexSet}(P_I)=\{D,E,F, G, H\}$.




\section{Algorithm}
\label{sec5}

We provide an pseudo-code of the algorithm~\ref{alg:NewIntegerHull} for {\verb|Integer Hull Computation|} using our procedure. The main steps of this algorithm rely on the function {\verb|ReplaceFacets|} as described in algorithm~\ref{alg:ReplaceFacets}.

\subsection{ReplaceFacets Algorithm}
In Line 1 and 2, we initialize the facets F to an empty list, and the integer vertices $V'$ to $V$. In Line 6 the {\verb|getCoeffs|} function calculates the equation of the line between vertices $va$ and $vb$. This forms the facets $F$ of $P$ in line 7.

In Line 10, we choose the vertex $v$ opposite to each facet $F$. In Line 12 - 16, we find the supporting hyperplane passing through $v$ and is parallel to $F$.

In Line 17 - 18, we find lines that are interior to the vertex before and after v, and contains integer points using our {\verb|nearestLine|} function. In Line 19 and 20, we use our {\verb|intersection|} function to find the point of intersection of these new lines with supporting hyperplane passing through $v$ and parallel to $F$.

In Line 21, we update the facet with F, and the intersection points to $V'$ in Line 22.

\subsection{IntegerHull Algorithm}
If the input polyhedral set $P$ is empty at Line 1, then we return an empty list in Line 2. Otherwise, in Line 3, we initialize the vertices and rays of $P$ to $V$ and $R$ respectively using Maple's {\verb|VerticesAndRays|} function, which is part of the {\verb|PolyhedralSets|} package.

If the number of vertices of $P$ is atleast $3$ then we implement our {\verb|SortPoints|}, {\verb|ReplaceFacets|} and {\verb|ReplaceVertices|} function, otherwise there can not be an integer hull and hence returns empty list in Line 10. The {\verb|SortPoints|} function sorts the vertices of $P$ based on angle with respect to a randomly selected origin vertex. The vertices obtained from the {\verb|ReplaceFacets|} might contain rational vertices, which {\verb|ReplaceVertices|} replaces with integer vertices that are contained within the polyhedral set. This is because the integer hull of a polyhedral set is contained within the polyhedral set itself.

In Line 11, we remove any duplicate vertices using Maple's {\verb|MakeUnique|} command.

In Line 12 if the number of integer vertices are more than $3$, we use Maple's {\verb|ConvexHull|} function that is part of the {\verb|ComputationalGeometry|} package, which computes the convex (integer) hull of the integer points. We return the integer hull of $P$ in Line 15. 

The integer hull of a polyhedral set is a polyhedral set. Therefore, the output of the algorithm is also a polyhedral set.

\begin{algorithm}[htb]
\small
\caption{{\sf ReplaceFacets}}\label{alg:ReplaceFacets}
\algnotext{EndIf}
\algnotext{EndFor}
\algnotext{EndWhile}
\begin{algorithmic}[1]
\Require {Vertices $V$ of $P$.}
\Ensure {Replace the facets that could not have integer point with the ones that could have.}

\State {{$F := []$}} \Comment{Initialize to empty list}
\State {{$V' := V$}} \Comment{Initialize to V}

\For {$i$ from $1$ to $|V|$}
    \State {$va := V'[i]$}
    \State {$vb := V'$[if$(i = |V|, 1, i + 1)$]}
    \State {$a, b, c$ := {\sf getCoeffs}($va, vb$)} \Comment{Coefficients of line $ax + cy = b$}
    \State {Append [$a, b, c$] to $F$} \Comment{Facets of $P$}
\EndFor

\For {$i$ from $1$ to $|V|$}
    \State {$a, b, c := F[i]$}
    \State {$val$ := if(ceil($|V|,2) + i > |V|$, modp(ceil($|V|/2) + i, |V|$), ceil($|V|/2) + i$)}
    \State {$v := V[val]$} \Comment{Choosing vertex which is opposite to the suppoting hyperplane}
    \If {$c = 0$}
        \State {$b := v[1]$} \Comment{Equation of a Parallel line through vertex $v$, when $ax = b$}
    \Else
        \State {$b' := v[2] + (a/c) \times v[1]$}
        \State {$b := c \times new_b$} \Comment{Equation of a Parallel line through vertex $v$}
    \EndIf
    \State {$a1, b1, c1$ := {\sf nearestLine}($[a, b, c], V$[if$(val = 1, |V|, val - 1)$])} \Comment{New line closer to $v$ with integer points}
    \State {$a2, b2, c2$ := {\sf nearestLine}($[a, b, c], V$[if$(val = |V|, 1, val + 1 )$])}

    \State {$p1$ := {\sf intersection}$([a, b, c], [a1, b1, c1])$} \Comment{Intersection of two lines}
    \State {$p2$ := {\sf intersection}$([a, b, c], [a2, b2, c2])$}

    \State {$F := [F[1], \ldots, F[i-1], [a, b, c], F[i+1], \ldots, F[|V|]$}
    \State {$V' := [V[1], \ldots, V[i-1], p1, p2, V'[i+2], \ldots, V[|V|]]$}
\EndFor

\State \Return {\sf $V'$} 

\end{algorithmic}
\end{algorithm}

\begin{algorithm}[htb]
\small
\caption{{\sf NewIntegerHull}}\label{alg:NewIntegerHull}
\algnotext{EndIf}
\algnotext{EndFor}
\algnotext{EndWhile}
\begin{algorithmic}[1]
\Require {$P$ a polyhedral set.}
\Ensure {Integer hull of the polyhedral set $P$.}

\If{$P=\phi$} 
    \State \Return {$[]$} \Comment{Returns empty set if empty input}
\EndIf
\State {{\sf V, R} := {\sf VerticesAndRays(P)}} \Comment{Initialize the vertices and rays of $P$}

\If{$\left|V\right| \geq 3$}
    \State {{\sf $\mathscr{V}$} := {\sf SortPoints(V, V[1])}} \Comment{Sorts the vertices of $P$ counter-clockwise starting with the first vertex} 
    \State {{\sf $\mathscr{V}$} := {\sf ReplaceFacets($\mathscr{V}$, P)}} 
    \State {{\sf $\mathscr{V}$} := {\sf ReplaceVertices($\mathscr{V}$)}} \Comment{Replace the vertices with integer vertices (inside of $P$)}
\Else
    \State {{\sf $\mathscr{V}$} := $[]$}
\EndIf

\State {{\sf $\mathscr{V}$}:= {\sf MakeUnique($\mathscr{V}$)}} \Comment{Remove the repeated integer vertices from the integer hull of $P$}

\If{$\left|\mathscr{V}\right| > 3$}
    \State {{\sf v} := {\sf ConvexHull($\mathscr{V}$)}} \Comment{Form the convex (integer) hull}
    \State { {\sf $\mathscr{V}$} := {\sf $\mathscr{V}$[v]}}
\EndIf




\State \Return {\sf PolyhedralSet($\mathscr{V}$, R)} \Comment{Returns the integer hull of $P$}

\end{algorithmic}
\end{algorithm}

\section{Experimentation}
\label{sec6}

In this section, we report on the software implementation of the algorithms proposed in the previous sections. We have implemented the algorithms in {\Maple} programming language. The {\Maple} version used is the 2024 release of {\Maple}. All the benchmarks are done on an Intel Core i7-7700T with Clockspeed: 2.9 GHz and Turbo Speed: 3.8 GHz. It has 4 cores and 8 threads.

\begin{tabular}{ |p{3cm}|p{3cm}|p{3cm}|p{3cm}|  }
 \hline
 \multicolumn{4}{|c|}{Preliminary Comparison Test} \\
 \hline
 Number of Vertices& Volume &New Algorithm&Existing Algorithm\\
 \hline
 10		&   29.9  				& 	663.00ms	&	823.00ms\\
 13		&	35.08				&	906.00ms	&	1.06ms\\
 13		&	40.56				&	996.00ms	&	1.10s\\
 12		&	40.63 				& 	898.00ms	&	928.00ms\\
 1000	&	69829.26 			&	855ms		& 	952ms\\
 15		&	263124.06 			&	1.65s 		& 	1.73s\\
 1000	& 	$6.54\times 10^6$ 	&	1.98s		& 	2.44s\\
 1500	&	$3.13\times 10^9$ 	& 	74.11s		&	89.72s\\
 \hline
\end{tabular}

\bibliographystyle{plain}
\bibliography{references}

\end{document}